\documentclass[11pt]{article}
\usepackage{latexsym, amsmath, amssymb, amsthm, a4, epsfig}
\usepackage{empheq}
\usepackage{graphicx}
\usepackage{color}
\usepackage{comment}
\usepackage{ifthen}
\usepackage{subcaption}
\usepackage{float}
\usepackage[
setpagesize=false,
bookmarksnumbered=true,%
bookmarksopen=true,%
colorlinks=true,%
linkcolor=blue,
citecolor=red,
]
{hyperref}

\newtheorem{theorem}{Theorem}[section]

\newtheorem{prop}[theorem]{Proposition}

\newtheorem{remark}{Remark}

\theoremstyle{definition}

\numberwithin{equation}{section}
\numberwithin{figure}{section}

\allowdisplaybreaks

\DeclareMathOperator{\disc}{disc}
\newcommand{\mua}{\mu_\text{a}}
\newcommand{\mut}{\mu_\text{t}}
\newcommand{\mus}{\mu_\text{s}}

\DeclareMathSymbol{\Real}{\mathalpha}{AMSb}{"52}

\begin{document}

\title{Remarks on propagation of discontinuities in stationary radiative transfer}

\author{Daisuke Kawagoe\thanks{Graduate School of Informatics, Kyoto University (d.kawagoe@acs.i.kyoto-u.ac.jp)}}


\date{\today}
\maketitle

\begin{abstract}
We consider the stationary transport equation with the incoming boundary condition. We are interested in discontinuities of the solution. Under the generalized convexity condition, it is known that it has only boundary-induced discontinuities, which are discontinuities arising from discontinuous boundary data, they propagate along positive characteristic lines, and we can reconstruct the attenuation coefficient from boundary measurements by the inverse X-ray transform. In this article, we observe that coefficient-induced discontinuities, discontinuities of the solution arising from discontinuous coefficients, would also appear without the generalized convexity condition. If the set of discontinuous points of the coefficients contains at most finite number of flat parts, coefficient-induced discontinuities do not affect the inverse X-ray transform. We also remark that, under the generalized convexity condition, a three dimensional inverse problem can be reduced to the two dimensional one. A numerical experiment is exhibited.
\end{abstract}

\section{Introduction} \label{sec:intro}

We consider the stationary transport equation
\begin{equation} \label{eq:STE}
\begin{split}
\xi \cdot \nabla_x I(x, \xi) + (\mua(x) + \mus(x)) I(x, \xi) = \mus(x) \int_{S^{d-1}} p(x, \xi, \xi^\prime) I(x, \xi^\prime)\,d\sigma_{\xi^\prime},\\ 
(x, \xi) \in \Omega \times S^{d-1}.
\end{split}
\end{equation}
It describes photon propagation in turbid media like biological tissue \cite{arridge,Arr}. The function $I(x, \xi)$ stands for density of photons at a point $x \in \Omega (\subset \mathbb{R}^d)$, $d = 2$ or $3$, with a direction $\xi \in S^{d-1}$. Here, $S^{d-1}$ is the unit sphere in $\mathbb{R}^d$. The coefficient $\mua$ characterizes absorption of particles by the media, and the coefficient $\mus$ and the integral kernel $p$ characterize scattering of particles in media; they are called the absorption coefficient, the scattering coefficient, and the scattering phase function, respectively. In this article, we let $\mut = \mua + \mus$ and call it the (total) attenuation coefficient. 

For a bounded domain $\Omega$ in $\mathbb{R}^d$ with $C^1$ boundary, define the incoming boundary $\Gamma_-$ and the outgoing boundary $\Gamma_+$ by
\begin{equation*}
\Gamma_{\pm} := \{ (x, \xi) \in \partial \Omega \times S^{d-1} \mid \pm n(x) \cdot \xi > 0 \},
\end{equation*}
where $n(x)$ is the outer unit normal vector at $x \in \partial \Omega$. The forward problem is to seek a solution $I$ to the equation \eqref{eq:STE} satisfying
\begin{equation} \label{eq:BC}
I(x, \xi) = I_0(x, \xi), \quad (x, \xi) \in \Gamma_-
\end{equation}
for a given function $I_0$ on $\Gamma_-$. 

We are interested in discontinuities of the solution $I$ to the boundary value problem \eqref{eq:STE}-\eqref{eq:BC}. In \cite{CK}, propagation of boundary-induced discontinuities, discontinuities of the solution arising from discontinuous boundary data, was discussed assuming that the coefficients (and the integral kernel) of the equation \eqref{eq:STE} were piecewise continuous (with respect to $x$) and that the collection of such pieces satisfies the generalized convexity condition \cite{1993Anik}. There, the authors made use of the propagation of boundary-induced discontinuities to obtain the X-ray transform of the attenuation coefficient $\mut$ from boundary measurements $(I_0, I|_{\Gamma_+})$. Since the coefficient $\mut$ is reconstructed by its X-ray transform, this study can be applied to solve an inverse problem to determine $\mut$ from the boundary measurements, which is a mathematical model of the optical tomography, a new medical imaging modality \cite{Arr}. 

The proposed idea in \cite{CK} was examined numerically \cite{CFK}. They exhibited some numerical experiments on the computerized tomography for reconstruction of $\mut$ in the unit disk. In one of the experiments, though $\mut$ was still piecewise continuous, the collection of the pieces violated the generalized convexity condition. Nevertheless, the numerical reconstruction seemed to be successful. This motivates us to discuss the case without the generalized convexity condition, or the set of discontinuous points of the coefficients has a flat part. We mention that the case where the set of the discontinuous points of the coefficients was a square was already discussed in \cite{CFK}.

The rest part of this article is organized as follows. In Section \ref{sec:BID}, we shall review results on propagation of boundary-induced discontinuities obtained in \cite{CK} as a preliminary. Thanks to the generalized convexity condition, only discontinuities of boundary data create those of the solution. If we remove the assumption of the generalized convexity condition, discontinuities of the coefficients would create those of the solution, which will be called coefficient-induced discontinuities in comparison to boundary-induced discontinuities. We shall describe coefficient-induced discontinuities in the two dimensional case in Section \ref{sec:CID_2D}, and extend the description to the three dimensional case in Section \ref{sec:CID_3D}. In addition, we shall show that, under the generalized convexity condition, the three dimensional inverse problem can be reduced to the two dimensional one, a tomographic imaging. We exhibit a numerical experiment for the computerized tomography at the end of Section \ref{sec:CID_3D}.

\section{Propagation of boundary-induced discontinuities} \label{sec:BID}

In this section, we review results on propagation of boundary-induced discontinuities obtained in \cite{CK}. 

\begin{figure}[h]
\begin{center}
\includegraphics[width=5cm]{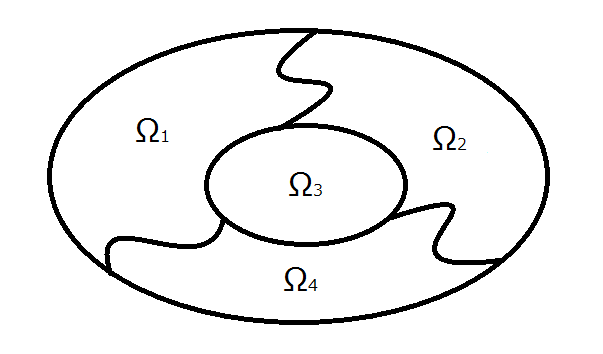}
\caption{An example of the generalized convexity condition}
\label{fig:gen_conv}
\end{center}
\end{figure}

We introduce a generalized convexity condition \cite{1993Anik}. Let $\Omega$ be a bounded convex domain in $\mathbb{R}^d$ with the $C^1$ boundary $\partial \Omega$. Suppose that $\overline{\Omega} = \cup_{j = 1}^N \overline{\Omega_j}$, where $\Omega_j$, $1 \leq j \leq N$, are disjoint subdomains of $\Omega$ with piecewise $C^1$ boundaries. Let $\Omega_0 := \cup_{j = 1}^N \Omega_j$. We say that the partition $\{ \Omega_j \}_{j = 1}^N$ satisfies the generalized convexity condition if, for all $(x, \xi) \in \Omega \times S^{d-1}$, the half line $\{ x - t\xi \mid t \geq 0 \}$ intersects with $\partial \Omega_0$ at most finite times. In other words, for all $(x, \xi) \in \Omega \times S^{d-1}$, there exist positive integer $l(x, \xi)$ and real numbers $\{ t_j(x, \xi) \}_{j = 1}^{l(x, \xi)}$ such that $0 \leq t_1(x, \xi) < t_2(x, \xi) < \cdots < t_{l(x, \xi)}(x, \xi)$, $x - t\xi \in \partial \Omega_0$ if and only if $t = t_j(x, \xi)$, and $\sup_{(x, \xi) \in \Omega \times S^{d-1}} l(x, \xi) < \infty$. In what follows, we use these notations $t_j(x, \xi)$ and $l(x, \xi)$ for the generalized convexity condition, and we put $t_0(x, \xi) = 0$. We put Figure \ref{fig:gen_conv} as an example of a partition of the domain satisfying the generalized convexity condition to facilitate readers' imagination. 

We next describe the setting in \cite{CK}. Let $\mut$ and $\mus$ be nonnegative bounded functions on $\Omega$ such that $\mut(x) \geq \mus(x)$ for all $x \in \Omega$. We assume that there exists a partition $\{ \Omega_j \}_{j=1}^N$ of the domain $\Omega$ such that all of discontinuous points of these two functions are contained in $\partial \Omega_0$ and it satisfies the generalized convexity condition. For the discussion in what follows, we reset $\mut(x) = \mus(x) = 0$ for $x \in \mathbb{R}^d \setminus \Omega_0$. It is worth mentioning that $\mut$ and $\mus$ are bounded continuous at least on $\Omega_0$.
We also assume that the integral kernel $p$ is a nonnegative bounded function on $\mathbb{R}^d \times S^{d-1} \times S^{d-1}$ which is continuous on $\Omega_0 \times S^{d-1} \times S^{d-1}$ and $p(x, \xi, \xi^\prime) = 0$ for $(x, \xi, \xi^\prime) \in (\mathbb{R}^d \setminus \Omega_0) \times S^{d-1} \times S^{d-1}$, and satisfies
\[
\int_{S^{d-1}} p(x, \xi, \xi^\prime)\,d\sigma_{\xi^\prime} = 1
\]
for all $(x, \xi) \in \Omega_0 \times S^{d-1}$. 

We regard the directional derivative $\xi \cdot \nabla_x I(x, \xi)$ as
\begin{equation*}
\xi \cdot \nabla_x I(x, \xi) := \left. \frac{d}{dt} I(x + t\xi, \xi) \right|_{t = 0}.
\end{equation*}
The measure $d\sigma_{\xi^\prime}$ is the Lebesgue measure on the unit sphere $S^{d-1}$.

We introduce some notations. Let
\begin{equation*}
D := (\Omega \times S^{d-1}) \cup \Gamma_-, \quad \overline{D} := D \cup \Gamma_+,
\end{equation*}
and we define two functions $\tau_\pm$ on $\overline{D}$ by
\begin{equation*}
\tau_{\pm}(x, \xi) := \inf \{ t > 0 \mid x \pm t\xi \not\in \Omega \}.
\end{equation*}
Let $\Gamma_{-, \xi}$ and $\Gamma_{-, x}$ be projections of $\Gamma_-$ on $\partial \Omega$ and $S^{d-1}$ respectively;
\begin{equation*}
\Gamma_{-, \xi} := \{ x \in \partial \Omega \mid n(x) \cdot \xi < 0 \}, \quad \xi \in S^{d-1},
\end{equation*}
and
\begin{equation*}
\Gamma_{-, x} := \{ \xi \in S^{d-1} \mid n(x) \cdot \xi < 0 \}, \quad x \in \partial \Omega.
\end{equation*}
Let $\disc(I)$ be the set of all discontinuous points for a function $I$.

We define a solution to the boundary value problem \eqref{eq:STE}-\eqref{eq:BC}. For all $(x, \xi) \in D$, integrating the equation \eqref{eq:STE} with respect to $x$ along the line $\{x - t\xi \mid t > 0\}$ until the line intersects with the boundary $\partial \Omega$ and taking the boundary condition \eqref{eq:BC} into consideration, we obtain the following integral equation:
\begin{equation} \label{eq:IE}
\begin{split}
I(x, \xi) =& \exp \Bigl(- M_t \bigl(x, \xi; \tau_-(x, \xi) \bigr) \Bigr) I_0(P(x, \xi), \xi)\\
&+ \int_0^{\tau_-(x, \xi)} \mus(x - s\xi) \exp \Bigl( - M_t (x, \xi; s) \Bigr)\\
& \quad \times \int_{S^{d-1}} p(x - s\xi, \xi, \xi^\prime) I(x - s\xi, \xi^\prime)\,d\sigma_{\xi^\prime}ds,
\end{split}
\end{equation}
where
\begin{equation} \label{def:Mt}
M_t (x, \xi; s) := \int_0^s \mut(x - r\xi)\,dr.
\end{equation}
We call a bounded measurable function $I$ on $D$ satisfying the integral equation \eqref{eq:IE} for all $(x, \xi) \in D$ a solution to the boundary value problem \eqref{eq:STE}-\eqref{eq:BC}. We remark that, under the generalized convexity condition on pieces of the coefficients, the solution satisfies the following properties: (i) it has the directional derivative $\xi \cdot \nabla_x I(x, \xi)$ at all $(x, \xi) \in \Omega_0 \times S^{d-1}$, (ii) it satisfies the stationary transport equation \eqref{eq:STE} for all $(x, \xi) \in \Omega_0 \times S^{d-1}$ and the boundary condition \eqref{eq:BC} for all $(x, \xi) \in \Gamma_-$, (iii) $I(\cdot, \xi)$ is continuous along the line $\{ x + t\xi \mid t \in \mathbb{R}\} \cap (\Omega \cup \Gamma_{-, \xi})$ for all $(x, \xi) \in D$, and (iv) $\xi \cdot \nabla_x I(\cdot, \xi)$ is continuous on the open line segments $\{ x - t\xi \mid t \in (t_{j-1}(x, \xi), t_j(x, \xi) ) \}$, $j = 1, \ldots, l(x, \xi)$ with $t_0(x, \xi) = 0$ for all $(x, \xi) \in \Omega_0 \times S^{d-1}$.

The first result shows that boundary-induced discontinuities propagate only along positive characteristic lines starting from discontinuous points of the incoming boundary data. Here, a positive characteristic line from a point $(x, \xi) \in \Gamma_-$ is defined by $\{ (x + t\xi, \xi) \mid t \geq 0\}$.

\begin{prop}[\cite{CK}] \label{prop:MR1}
Suppose that a boundary data $I_0$ is bounded and that it satisfies at least one of the following two conditions.
\begin{enumerate}
\item $I_0(x, \cdot)$ is continuous on $\Gamma_{-, x}$ for almost all $x \in \partial \Omega$,
\item $I_0(\cdot, \xi)$ is continuous on $\Gamma_{-, \xi}$ for almost all $\xi \in S^{d-1}$.
\end{enumerate}
Then, there exists a unique solution $I$ to the boundary value problem $\eqref{eq:STE}$-$\eqref{eq:BC}$, and we have
\begin{equation*}
\disc(I) = \{ (x_* + t\xi_*, \xi_*) \mid (x_*, \xi_*) \in \disc(I_0), 0 \leq t < \tau_+(x_*, \xi_*) \}.
\end{equation*}
\end{prop}

\begin{remark}
Proposition \ref{prop:MR1} implies that, for a bounded continuous boundary data $I_0$ on $\Gamma_-$, the unique solution $I$ is bounded continuous on $D$.
\end{remark}

\begin{remark}
Anikonov et al.\ \cite{1993Anik} showed Proposition \ref{prop:MR1} with the condition 2. Our main contribution is to show Proposition \ref{prop:MR1} with the condition 1.
\end{remark}

As the second result, we shall discuss boundary-induced discontinuities of the solution extended up to $\Gamma_+$. In other words, we can extend the domain of the solution $I$ up to $\Gamma_+$, and we see that boundary-induced discontinuities propagate along positive characteristic lines up to $\Gamma_+$.

\begin{prop}[\cite{CK}] \label{prop:MR2}
Let a boundary data $I_0$ satisfy assumptions in Proposition \ref{prop:MR1} and let $I$ be the solution to the boundary value problem $\eqref{eq:STE}$-$\eqref{eq:BC}$. Then, it can be extended up to $\Gamma_+$, which is denoted by $\overline{I}$, by
\begin{equation*}
\overline{I}(x, \xi) :=
\begin{cases}
I(x, \xi), \quad &(x, \xi) \in D, \\
\displaystyle \lim_{t \downarrow 0} I(x - t\xi, \xi), \quad &(x, \xi) \in \Gamma_+.
\end{cases}
\end{equation*}
Moreover, we have
\begin{equation*}
\disc(\overline{I}) = \{ (x_* + t\xi_*, \xi_*) \mid (x_*, \xi_*) \in \disc(I_0), 0 \leq t \leq \tau_+(x_*, \xi_*) \}.
\end{equation*}
\end{prop}

We state the decay of boundary-induced discontinuities in some situation as the third result. Let $\gamma$ be two points in $\partial \Omega$ when $d = 2$, while let $\gamma$ be a simple closed curve in $\partial \Omega$ when $d = 3$. Then, $\gamma$ splits $\partial \Omega$ into two connected components $A$ and $B$, that is $\partial \Omega = A \cup B \cup \gamma$ and $A \cap B = A \cap \gamma = B \cap \gamma = \varnothing$. 
We put an incoming boundary data $I_0$ by
\begin{equation} \label{eq:BCJ}
I_0(x, \xi) =
\begin{cases}
C, \quad (x, \xi) \in ( (A \cup \gamma) \times S^{d-1} ) \cap \Gamma_-,\\
0, \quad (x, \xi) \in (B \times S^{d-1} ) \cap \Gamma_-,
\end{cases}
\end{equation}
where $C$ is a non-zero constant. We note that $I_0$ satisfies the condition 1 of Proposition \ref{prop:MR1}, and that $\disc(I_0) = \{ (x_*, \xi_*) \mid x_* \in \gamma, \xi_* \in \Gamma_{-, x_*} \}$.

For $(\overline{x}, \overline{\xi}) \in \disc(\overline{I})$, we define a jump $[\overline{I}] (\overline{x}, \overline{\xi})$ by
\[
[\overline{I}](\overline{x}, \overline{\xi}) := \lim_{\substack{x \rightarrow \overline{x} \\ P(x, \overline{\xi}) \in (A \cup \gamma)}} I(x, \overline{\xi}) - \lim_{\substack{x \rightarrow \overline{x} \\ P(x, \overline{\xi}) \in B}} I(x, \overline{\xi}),
\]
where
\begin{equation*}
P(x, \xi) := x - \tau_-(x, \xi)\xi.
\end{equation*}
We note that, in our situation, $[I_0](x, \xi) = C$ for all $(x, \xi) \in \disc(I_0) = (\gamma \times S^{d-1}) \cap \Gamma_-$. In this situation, we have the following theorem.

\begin{prop}[\cite{CK}] \label{prop:MR3}
Let $\overline{I}$ be the extended solution to the boundary value problem \eqref{eq:STE}-\eqref{eq:BC} with the incoming boundary data given by \eqref{eq:BCJ}, and let $(x^*, \xi^*) \in \disc(\overline{I})$. Then, we have
\begin{equation} \label{eq:decay_jump}
[\overline{I}](x^*, \xi^*) = C \exp \left( - \int_0^{\tau_-(x^*, \xi^*)} \mut(x^* - r\xi^*)\,dr \right).
\end{equation}
\end{prop}

In particular, we take a point $(x^*, \xi^*) \in \disc(\overline{I}) \cap \Gamma_+$. By Proposition \ref{prop:MR3}, we have
\begin{equation} \label{eq:Xray}
X \mut (x^*, \xi^*) = \int_0^{\tau_-(x^*, \xi^*)} \mut(x^* - r\xi^*)\,dr = - \log \left( [\overline{I}](x^*, \xi^*) / C \right).
\end{equation} 
The right hand side is obtained from observed data. By changing a position of $\gamma$, we can observe the X-ray transform of $\mut$. Then, applying a well-known method such as filtered back projection \cite{Nat}, we can reconstruct the attenuation coefficient $\mut$.

We give a sketch of proofs of the above proposition for readers' convenience. Define a sequence of functions $\{I^{(n)}\}_{n \geq 0}$ on $D$ by
\begin{equation} \label{eq:F0}
I^{(0)}(x, \xi) := \exp \Bigl(- M_t(x, \xi; \tau_- (x, \xi)) \Bigr) I_0(P(x, \xi), \xi) 
\end{equation}
and
\begin{equation} \label{eq:F1}
\begin{split}
I^{(n+1)}(x, \xi) :=& \int_0^{\tau_-(x, \xi)} \mus(x - s\xi) \exp \Bigl(- M_t (x, \xi; s) \Bigr) \\
&\quad \times \int_{S^{d-1}} p(x - s\xi, \xi, \xi^\prime) I^{(n)}(x - s\xi, \xi^\prime)\,d\sigma_{\xi^\prime}ds.
\end{split}
\end{equation}
Then, we have
\[
\sup_{(x, \xi) \in D} |I^{(0)}(x, \xi)| = \sup_{(x, \xi) \in \Gamma_-} |I_0(x, \xi)|
\]
and
\[
\sup_{(x, \xi) \in D} |I^{(n + 1)}(x, \xi)| \leq M \sup_{(x, \xi) \in D} |I^{(n)}(x, \xi)|,
\]
where 
\[
M := \sup_{(x, \xi) \in D} \biggl( 1 - \exp \Bigl(-M_t \bigl(x, \xi; \tau_-(x, \xi) \bigr) \Bigr) \biggr).
\]
Since $M < 1$, the sum $I := \sum_{n = 0}^\infty I^{(n)}$ converges absolutely and uniformly in $D$, and it is the unique solution to the integral equation \eqref{eq:IE}. 

We decompose the solution $I$ to the integral equation \eqref{eq:IE} into two parts as the following:
\[
I(x, \xi) = I_D(x, \xi) + I_C(x, \xi),
\]
where
\begin{align*}
I_D(x, \xi) :=& \exp \Bigl(- M_t(x, \xi; \tau_- (x, \xi)) \Bigr) I_0(P(x, \xi), \xi),\\ 
I_C(x, \xi) :=& \sum_{n = 1}^\infty I^{(n)}(x, \xi). 
\end{align*}
We show that the discontinuous points of the function $I_D$ is described as
\[
\disc(I_D) = \{(x_* + t\xi_*, \xi_*) \mid (x_*, \xi_*) \in \disc(I_0), 0 \leq t < \tau_+(x_*, \xi_*) \},
\]
while the function $I_C$ is bounded continuous on $D$. 

Under the assumption of the generalized convexity condition, the function $M_t(\cdot, \cdot; s)$ is continuous on $D$ for all $s \in [0, R]$, where $R$ is the diameter of the domain $\Omega$. In particular, since $\mut = 0$ on $\mathbb{R}^d \setminus \Omega_0$, it follows that the function $M_t(\cdot, \cdot; \tau_-(\cdot, \cdot))$ is continuous on $D$. Thus, discontinuities of the function $I_D$ only come from those of the boundary data $I_0$.

We discuss continuity of the function $I^{(1)}$. Let $G$ be the function defined by
\begin{equation} \label{eq:G}
G(x, \xi) := \int_{S^{d-1}} p(x, \xi, \xi^\prime) I_D(x, \xi^\prime)\,d\sigma_{\xi^\prime}.
\end{equation}
More explicitly, we have
\begin{equation} \label{eq:formula_G}
\begin{split}
G(x, \xi) =& \int_{S^{d-1}} p(x, \xi, \xi^\prime) \exp \Bigl(- M_t (x, \xi^\prime; \tau_- (x, \xi^\prime)) \Bigr) I_0(P(x, \xi^\prime), \xi^\prime) \,d\sigma_{\xi^\prime}\\
=& \int_{\partial \Omega} p \left(x, \xi, \frac{x - y}{|x - y|} \right) \exp \left(- M_t \left(x, \frac{x - y}{|x - y|}; |x - y|\right) \right)\\
& \times I_0 \left(y, \frac{x - y}{|x - y|} \right) \dfrac{ |n(y) \cdot (x - y) |} {|x - y|^d} \,d\sigma_y.
\end{split}
\end{equation}
Although the function $I_0$ may be discontinuous, under the assumption of Proposition \ref{prop:MR1}, the function $G$ is bounded continuous on $\Omega_0 \times S^{d-1}$. Since the function $I^{(1)}$ is described as
\[
I^{(1)}(x, \xi) = \int_0^{\tau_-(x, \xi)} \mus(x - s\xi) \exp \Bigl(- M_t (x, \xi; s) \Bigr) G(x - s\xi, \xi)\,ds,
\]
it is bounded continuous on $D$ by the dominated convergence theorem. We inductively say that the functions $\{ I^{(n)} \}_{n \geq 1}$ are bounded continuous on $D$. Since the series $\sum_{n = 1}^\infty I^{(n)}$ converges uniformly on $D$, the function $I_C$ is also bounded continuous on $D$.

We consider the extension of functions $I_D$ and $I_C$ up to $\Gamma_+$ in the way stated in Proposition \ref{prop:MR2}, which are denoted by $\overline{I_D}$ and $\overline{I_C}$ respectively. The function $\overline{I_C}$ is bounded continuous on $\overline{D}$ by the dominated convergence theorem, while the boundary-induced discontinuity in $I_D$ propagates up to $\Gamma_+$. 

Now we discuss the jump of the extended solution in Proposition \ref{prop:MR3}. Let $(x^*, \xi^*) \in \disc(\overline{I})$. By Proposition \ref{prop:MR2}, we have $(P(x^*, \xi^*), \xi^*) \in \disc(I_0)$ or $P(x^*, \xi^*) \in \gamma$. Thus, we have
\begin{align*}
[\overline{I_D}](x^*, \xi^*) =& \exp \left( - M_t(x^*, \xi^*; \tau_-(x^*, \xi^*)) \right)[I_0](P(x^*, \xi^*), \xi^*)\\ 
=& C \exp \left( - M_t(x^*, \xi^*; \tau_-(x^*, \xi^*)) \right).
\end{align*}
On the other hand, since $\overline{I_C}$ is bounded continuous on $\overline{D}$, we have $[\overline{I_C}](x^*, \xi^*) = 0$. Thus, we have
\begin{align*}
[\overline{I}](x^*, \xi^*) =& [\overline{I_D}](x^*, \xi^*) + [\overline{I_C}](x^*, \xi^*)\\ 
=& C \exp \left( - \int_0^{\tau_-(x^*, \xi^*)} \mut(x^* - r\xi^*)\,dr \right).
\end{align*}

This is the sketch of proofs from Proposition \ref{prop:MR1} to Proposition \ref{prop:MR3}.

\section{Coefficient-induced discontinuities in 2D} \label{sec:CID_2D}

In Section \ref{sec:BID}, we reviewed propagation of boundary-induced discontinuities assuming that the set of discontinuous points of the coefficients was contained in the boundary of a partition $\{ \Omega_j \}_{j = 1}^N$ of the domain $\Omega$ and the partition satisfied the generalized convexity condition. In this section, we investigate discontinuities of solutions without the generalized convexity condition in the two dimensional case.

We introduce a setting for this section. Let $\Omega$ be a bounded convex domain in $\mathbb{R}^2$ with $C^1$ boundary. Suppose that $\mut$, $\mus$ and $p$ satisfy the assumptions in Section \ref{sec:BID} except that the partition $\{ \Omega_j \}_{j = 1}^N$ satisfies the generalized convexity condition. Instead, we assume that $\partial \Omega_0 \setminus \partial \Omega$ contains at most finite number of (open) line segments $\{ L_j \}$, $j = 1, \ldots, K$. Since the case $K = 0$ corresponds to the generalized convexity condition, we only consider the case $K \geq 1$. We note that this assumption holds if $\partial \Omega_0$ is piecewise $C^1$ and it contains a line segment.

For $1 \leq j \leq K$, let $\xi^{(j)}_+ \in S^1$ be a direction of the line segment $L_j$. Namely, for any $x \in L_j$, there exists $t > 0$ such that $x + t \xi^{(j)}_+ \in L_j$. We also let $\xi^{(j)}_- := - \xi^{(j)}_+$. Then, Proposition \ref{prop:MR1} is extended as follows.

\begin{theorem} \label{theorem:MR1}
Suppose that a boundary data $I_0$ is bounded and that it satisfies at least one of the following two conditions.
\begin{enumerate}
\item $I_0(x, \cdot)$ is continuous on $\Gamma_{-, x}$ for almost all $x \in \partial \Omega$,
\item $I_0(\cdot, \xi)$ is continuous on $\Gamma_{-, \xi}$ for almost all $\xi \in S^1$.
\end{enumerate}
Then, there exists a unique solution $I$ to the boundary value problem $\eqref{eq:STE}$-$\eqref{eq:BC}$. Moreover, we have
\[
\disc(I) \subset \disc_B \cup \disc_C,
\]
where
\[
\disc_B := \{ (x_* + t\xi_*, \xi_*) \mid (x_*, \xi_*) \in \disc(I_0), 0 \leq t < \tau_+(x_*, \xi_*) \}
\]
and
\[
\disc_C := \cup_{j = 1}^K \{ (x + t\xi^{(j)}_\pm, \xi^{(j)}_\pm) \mid x \in L_j, 0 \leq t < \tau_+(x, \xi^{(j)}_\pm) \}.
\]
\end{theorem}

We call $\disc_B$ boundary-induced discontinuities as we saw in Section \ref{sec:BID} since it is characterized by $\disc(I_0)$. In contrast, the set $\disc_C$ is characterized by $\{ L_j \}$ in $\partial \Omega_0$, which are related to discontinuous points of the coefficients. Thus, we call it coefficient-induced discontinuities. 

\begin{proof}
Even without the generalized convexity condition, the sum $I = \sum_{n = 0}^\infty I^{(n)}$ with functions $\{ I^{(n)} \}$ defined by \eqref{eq:F0}-\eqref{eq:F1} converges absolutely and uniformly in D, and it is the unique solution to the boundary value problem \eqref{eq:STE}-\eqref{eq:BC}. In what follows, we discuss its discontinuities.

We first discuss discontinuities of the function $I^{(0)}$, which consists of the boundary data $I_0$ and the exponential factor with the function $M_t$. As we saw in Section \ref{sec:BID}, the discontinuous boundary data $I_0$ creates discontinuities of $I^{(0)}$. We focus on discontinuities of the function $M_t$.

By the dominated convergence theorem, the function $M_t$ defined by \eqref{def:Mt} is continuous at $(x, \xi, s) \in \Omega \times S^1 \times [0, R]$ if $\xi \neq \xi^{(j)}_\pm$, $j = 1, \ldots, K$. On the other hand, $M_t(\cdot, \xi^{(j)}_\pm; s)$ may be discontinuous on positive characteristic lines
\[
\{ x + t \xi^{(j)}_\pm \mid x \in L_j, 0 \leq t < \tau_-(x, \xi^{(j)}_\pm) \}
\]
for some large $s \in [0, R]$. Combining the above cases, we obtain $\disc(I^{(0)}) \subset \disc_B \cup \disc_C$. Here, the equality ``$=$'' does not hold in general due to three reasons. First, we assumed that all of discontinuous points of the coefficients were included in $\partial \Omega_0$, which does not imply that all of the point $\partial \Omega_0$ are discontinuous points of the coefficients. Second, a coefficient-induced discontinuity does not appear if the boundary data on the corresponding characteristic line is zero. Third, cancellations of boundary-induced discontinuities and coefficient-induced discontinuities could occur.

We next investigate discontinuities of the other functions $I^{(n)}$. For $n \geq 0$, let
\[
G^{(n)}(x, \xi) := \int_{S^1} p(x, \xi, \xi') I^{(n)}(x, \xi')\,d\sigma_{\xi'}
\]
so that 
\[
I^{(n + 1)}(x, \xi) = \int_0^{\tau_-(x, \xi)} \mus(x - s\xi) \exp \Bigl( -M_t(x, \xi; s) \Bigr) G^{(n)}(x - s\xi, \xi)\,ds.
\]
We remark that the function $G^{(0)}$ is the function $G$ defined by \eqref{eq:G}. Thanks to the formula \eqref{eq:formula_G} and the dominated convergence theorem, we can show that $G^{(0)}$ is bounded continuous on $\Omega_0 \times S^1$. Thus, discontinuities of the function $I^{(1)}$ only come from those of coefficients. In other words, we have $\disc(I^{(1)}) \subset \disc_C$. Inductively, by applying the dominated convergence theorem again, we see that $G^{(n)}$ is bounded continuous on $\Omega_0 \times S^1$ for all $n \geq 1$ and $\disc(I^{(n)}) \subset \disc_C$ for all $n \geq 2$.

Since all of functions $\{ I^{(n)} \}$ are bounded continuous at least on $D \setminus (\disc_B \cup \disc_C)$ and the sum $I = \sum_{n = 0}^\infty I^{(n)}$ converges absolutely and uniformly in $D$, it is also bounded continuous at least on $D \setminus (\disc_B \cup \disc_C)$, which means that
\[
\disc(I) \subset \disc_B \cup \disc_C.
\]
This completes the proof.
\end{proof}

The following theorem can be shown in the same way as in the proof of Proposition \ref{prop:MR2}.

\begin{theorem} \label{theorem:MR2}
Let a boundary data $I_0$ satisfy assumptions in Theorem \ref{theorem:MR1} and let $I$ be the solution to the boundary value problem $\eqref{eq:STE}$-$\eqref{eq:BC}$. Then, it can be extended up to $\Gamma_+$, which is denoted by $\overline{I}$, by
\begin{equation*}
\overline{I}(x, \xi) :=
\begin{cases}
I(x, \xi), \quad &(x, \xi) \in D, \\
\displaystyle \lim_{t \downarrow 0} I(x - t\xi, \xi), \quad &(x, \xi) \in \Gamma_+.
\end{cases}
\end{equation*}
Moreover, we have
\begin{equation*}
\disc(\overline{I}) \subset \overline{\disc_B} \cup \overline{\disc_C},
\end{equation*}
where
\[
\overline{\disc_B} := \{ (x_* + t\xi_*, \xi_*) \mid (x_*, \xi_*) \in \disc(I_0), 0 \leq t \leq \tau_+(x_*, \xi_*) \}
\]
and
\[
\overline{\disc_C} := \bigcup_{j = 1}^K \{ (x + t\xi^{(j)}_\pm, \xi^{(j)}_\pm) \mid x \in L_j, 0 \leq t \leq \tau_+(x, \xi^{(j)}_\pm) \}.
\]
\end{theorem}

We can also discuss jumps of boundary-induced discontinuities with the boundary condition \eqref{eq:BCJ}.

\begin{theorem} \label{theorem:MR3}
Let $\overline{I}$ be the extended solution to the boundary value problem \eqref{eq:STE}-\eqref{eq:BC} with the incoming boundary data given by \eqref{eq:BCJ}, and let $(x^*, \xi^*) \in \disc(\overline{I})$. Suppose that $\xi^* \neq \xi^{(j)}_\pm$, $j = 1, \ldots, K$. Then, the identity \eqref{eq:decay_jump} holds.
\end{theorem}

We remark that we can still reconstruct the coefficient $\mut$ from observed data. By Theorem \ref{theorem:MR3}, for $(x^*, \xi^*) \in \disc(\overline{I}) \cap \Gamma_+$ with $\xi^* \neq \xi^{(j)}_\pm$, $j = 1, \ldots, K$, we can observe the X-ray transform of $\mut$ \eqref{eq:Xray}. Since the set of the irregular directions $\{ \xi^{(j)}_\pm \}_{j = 1, \ldots, K}$ is a null set in $S^1$, lack of the X-ray transform on the set (or on $(x^*, \xi^{(j)}_\pm) \in \disc(\overline{I}) \cap \Gamma_+$ more precisely) does not affect the inverse X-ray transform. Thus, we can apply the inverse X-ray transform to the observation \eqref{eq:Xray} to reconstruct the attenuation coefficient $\mut$. 

\section{Remarks on discontinuities in 3D} \label{sec:CID_3D}

In this section, we make two remarks on propagation of discontinuities in the three dimensional case.

The first one is on counterparts of results in Section \ref{sec:CID_2D}. Let $\Omega$ be a bounded convex domain in $\mathbb{R}^3$ with $C^1$ boundary. Suppose that there exists a partition $\{ \Omega_j \}_{j = 1}^N$ of the domain $\Omega$ with piecewise $C^1$ boundaries such that $\mut$, $\mus$ and $p$ are bounded continuous on each $\Omega_j$, $\partial \Omega_0 \setminus \partial \Omega$ contains at most finite number of (closed) flat parts $\{ F_j \}$, $j = 1, \ldots, K$ and that the other parts $\partial \Omega_0 \setminus (\partial \Omega \cup_{j = 1}^K F_j)$ do not contain any line segments. In this case, the same results as in Section \ref{sec:CID_2D} hold. 

For $1 \leq j \leq K$, let $n_j$ be a normal vector to $F_j$ and let $S^2_j := \{ \xi \in S^2 \mid \xi \cdot n_j = 0 \}$. Then, we have
\[
\disc(I) \subset \disc_B \cup \disc_{C, 3d},
\]
where
\[
\disc_{C, 3d} := \bigcup_{j = 1}^K \{ (x + t\xi, \xi) \mid x \in F_j, \xi \in S^2_j, 0 \leq t < \tau_+(x, \xi) \}.
\]
Also, for the extended solution $\overline{I}$, we have
\[
\disc(\overline{I}) \subset \overline{\disc_B} \cup \overline{\disc_{C, 3d}},
\]
where
\[
\overline{\disc_{C, 3d}} := \bigcup_{j = 1}^K \{ (x + t\xi, \xi) \mid x \in \overline{F_j}, \xi \in S^2_j, 0 \leq t \leq \tau_+(x, \xi) \}.
\]
Moreover, with the boundary condition \eqref{eq:BCJ}, the identity \eqref{eq:decay_jump} holds for all $(x^*, \xi^*) \in \disc(\overline{I})$ with $\xi^* \notin \cup_{j = 1}^K S^2_j$. Since the set $\cup_{j = 1}^K S^2_j$ is a null set in $S^2$, we can apply the inverse X-ray transform to observed data to reconstruct the coefficient $\mut$.

The second one is on the reduction of the three dimensional inverse problem to the two dimensional one. For the sake of simplicity, we assume that the partition $\{ \Omega_j \}_{j = 1}^N$ satisfies the generalized convexity condition. Since the domain $\Omega$ is bounded, there exist two constants $a, b \in \mathbb{R}$ such that it is contained in the slab $\mathbb{R}^2 \times (a, b)$. For a fixed $x_0 \in (a, b)$, take
\[
A := \partial \Omega \cap \{ x_3 > x_0 \}, \quad B := \partial \Omega \cap \{ x_3 < x_0 \}, \quad \gamma := \partial \Omega \cap \{ x_3 = x_0 \},
\]
and consider the boundary condition \eqref{eq:BCJ}. Furthermore, let $\tilde{\Omega} := \Omega \cap \{ x_3 = x_0 \}$ and $\tilde{\Omega}_j := \Omega_j \cap \{ x_3 = x_0 \}$, $j = 1, \ldots, N$. We remark that $\{ \tilde{\Omega}_j \}_{j = 1}^N$ is a partition of $\tilde{\Omega}$ and it also satisfies the generalized convexity condition.

Let $S^2_0 := \{ (\xi_1, \xi_2, \xi_3) \in S^2 \mid \xi_3 = 0 \}$. Then, by Proposition \ref{prop:MR2}, for each $\xi \in S^2_0$, we have $\disc(\overline{I}(\cdot, \xi)) = \overline{\tilde{\Omega}}$. In other words, if we restrict the direction $\xi$ to $S^2_0$, boundary-induced discontinuities propagate in $\tilde{\Omega}$. Moreover, we can observe the jump of the extended solution $\overline{I}$ \eqref{eq:decay_jump} on $\Gamma_+ \cap (\gamma \times S^2_0)$. Since the X-ray transform of the coefficient $\mut$ on $\{ x_3 = x_0 \} \times S^2_0$ is obtained, we can reconstruct the restriction of the coefficient $\mut|_{\tilde{\Omega}}$ from boundary measurements. It is worth mentioning that we do not need to rearrange the closed curve $\gamma$ in this reconstruction procedure.

\begin{figure}[h]
\begin{center}
\includegraphics[width=5cm]{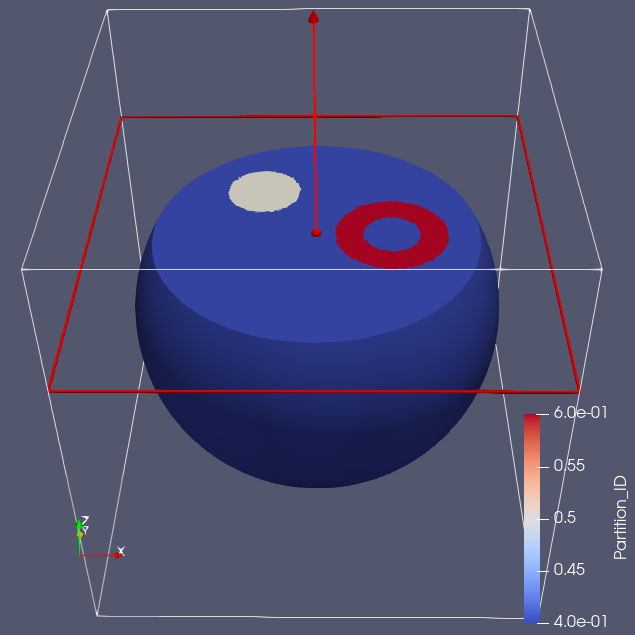}
\caption{The section of the example at $x_3 = 0.50$}
\label{fig:sec_0.5}
\end{center}
\end{figure}

We exhibit a numerical experiment for the second remark. Let $\Omega = \{ x \in \mathbb{R}^3 \mid |x| < 1 \}$ be the unit ball. Also, let 
\begin{align*}
\Omega_1 :=& \{ (x_1, x_2, x_3) \in \Omega \mid 0.031 < (x_1 - 0.43)^2 + x_2^2 + (x_3 - 0.50)^2 < 0.13 \},\\
\Omega_2 :=& \{ (x_1, x_2, x_3) \in \Omega \mid (x_1 + 0.38)^2 + (x_2 - 0.38)^2 + (x_3 - 0.50)^2 < 0.047 \},\\
\Omega_3 :=& \Omega \setminus (\overline{\Omega_1} \cup \overline{\Omega_2}).
\end{align*}
The section of the domain at $x_3 = 0.50$ with the partition is shown as Figure \ref{fig:sec_0.5}. We remark that $\Omega_1$, $\Omega_2$ and $\Omega_3$ correspond to the red, white and blue parts in Figure \ref{fig:sec_0.5} respectively.

We set $\mus = 0.3$ on $\Omega$ and define $\mua$ by
\[
\mua(x) = 
\begin{cases}
0.3, &x \in \Omega_1,\\
0.2, &x \in \Omega_2,\\
0.1, &x \in \Omega_3.
\end{cases}
\]
We adopt the scattering phase function $p$ of the following form:
\[
p(x, \xi, \xi') := \frac{1 - g^2}{4 \pi} \frac{1}{(1 + g^2 - 2g \cos(\xi \cdot \xi'))^{3/2}}
\]
with $g = 0.9$. The function is called the Henyey-Greenstein kernel, and it is typically used in the field of the optical tomography \cite{HG}.

\begin{figure}
\begin{center}
\includegraphics[width=0.9\columnwidth]{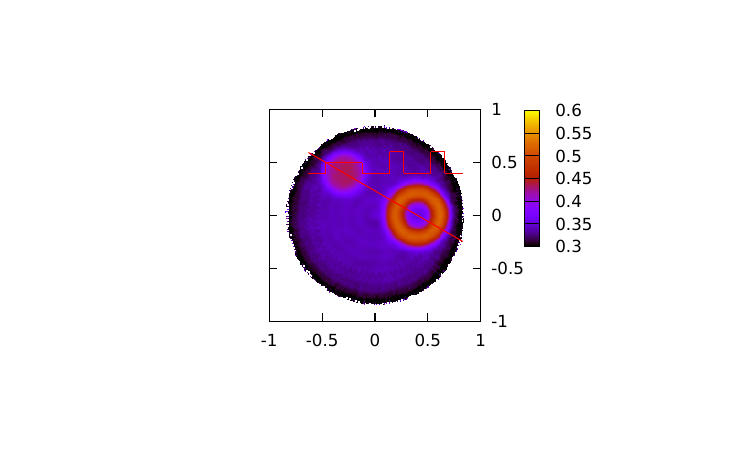}
    \caption{The reconstructed image of $\mut$}
    \label{fig:a}
\end{center}
\end{figure}

The reconstructed image of $\mut (=\mus + \mua)$ based on the proposed method is shown in Figure \ref{fig:a}. The A-analytic theory was also applied to stabilize the reconstruction process \cite{FST}.

\begin{figure}
\begin{center}
\includegraphics[width=0.8\columnwidth]{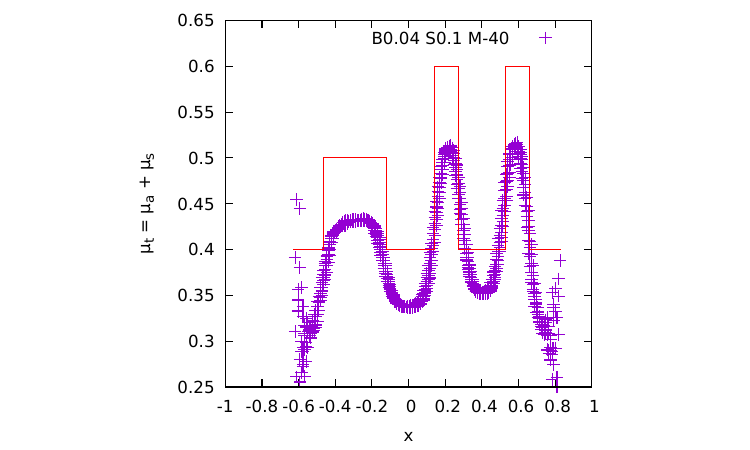}
    \caption{The section of the reconstructed image}
    \label{fig:b}
\end{center}
\end{figure}

Figure \ref{fig:b} shows the section of the reconstructed image along a red diagonal line in Figure \ref{fig:a}. 

\section*{Acknowledgement}
This work was supported by JSPS KAKENHI grant number JP24K00539. The author would like to thank Professor Hiroshi Fujiwara for his numerical reconstruction based on the propagation of discontinuities and for providing figures from Figure \ref{fig:sec_0.5} to Figure \ref{fig:b} in this article.


\begin{thebibliography}{11}

\bibitem{1993Anik}
D.~S.~Anikonov, I.~V.~Prokhorov and A.~E.~Kovtanyuk, 
Investigation of scattering and absorbing media by the methods of X-ray tomography, 
J. Inv. Ill-Posed Problems, 
\textbf{1}, 
(1993),
259--281.

\bibitem{arridge}
S.~R.~Arridge, 
Optical tomography in medical imaging,
Inverse Problems, 
\textbf{15}, 
(1999),
41--93.

\bibitem{Arr}
S.~R.~Arridge and J.~C.~Schotland,
Optical tomography: Forward and inverse problems,
Inverse Problems, 
\textbf{25}, 
(2009), 
123010, 
59pp.

\bibitem{CFK}
I.~Chen, H.~Fujiwara and D.~Kawagoe,
Tomography from scattered signals obeying the stationary radiative transport equation,
Practical Inverse Problems and Their Prospects,
(2023),
27--46.

\bibitem{CK} 
I.~Chen and D.~Kawagoe, 
Propagation of boundary-induced discontinuity in stationary radiative transfer and its application to the optical tomography, 
Inverse Probl. Imaging, 
\textbf{13}(2), 
(2019),
337--351. 

\bibitem{FST}
H.~Fujiwara, K.~Sadiq and A.~Tamasan,
A Fourier approach to the inverse source problem in an absorbing and anisotropic scattering medium,
Inverse Problems,
\textbf{36},
(2020),
33 pp.

\bibitem{HG}
L.~G.~Henyey and J.~L.~Greenstein,
Diffuse radiation in the galaxy,
Astrophys. J.,
\textbf{93},
(1941),
70--83.

\bibitem{Nat}
F.~Natterer,
The Mathematics of Computerized Tomography,
SIAM, 
Germany,
2001.

\end{thebibliography}
\end{document}